\documentclass[a4paper,12pt,reqno]{amsart}
\usepackage[a4paper,top=2cm,bottom=2cm,left=2.2cm,right=2.2cm]{geometry}

\usepackage{orcidlink}
\usepackage{microtype}
\usepackage{mlmodern}

\usepackage[utf8]{inputenc}
\usepackage[english]{babel}

\usepackage{amsmath,amssymb,amsfonts,amsthm}
\usepackage{graphicx,tikz}
\usetikzlibrary{arrows,matrix}
\usepackage[table]{xcolor}
\usepackage[centertableaux]{ytableau}
\usepackage{array}
\usepackage{booktabs}
\usepackage{tabularx}
\usepackage{hyperref}

\usepackage{longtable}
\usepackage{makecell}

\theoremstyle{plain}
\newtheorem{theorem}{Theorem}[section]

\newtheorem{lemma}[theorem]{Lemma}
\newtheorem{conjecture}[theorem]{Conjecture}
\newtheorem{proposition}[theorem]{Proposition}

\theoremstyle{definition}

\theoremstyle{remark}
\newtheorem{remark}[theorem]{Remark}

\numberwithin{equation}{section}

\title[Partitions with Parts Separated by Parity]{Inequalities and Parity Biases involving Partitions with Parts Separated by Parity}

\author[A. M. Alanazi]{Abdulaziz M. Alanazi\,\orcidlink{0000-0002-7909-2539}}
\address[A. M. Alanazi]{Department of Mathematics, Faculty of Sciences, University of Tabuk, P.O.Box 741, Tabuk 71491, Saudi Arabia}
\email{am.alenezi@ut.edu.sa}

\author[G. Kalita]{Gitartha Kalita\, \orcidlink{0009-0004-6073-9584}}
\address[G. Kalita]{Department of Mathematical Sciences, Tezpur University, Napaam, 784028, Assam, India}
\email{gitartha.kalita062000@gmail.com}

\author[P. J. Mahanta]{Pankaj Jyoti Mahanta\,\orcidlink{0000-0001-6460-7395}}
\address[P. J. Mahanta]{Department of Mathematical Sciences, Tezpur University, Napaam, 784028, Assam, India}
\email{pjm2099@gmail.com}

\author[M. P. Saikia]{Manjil P. Saikia\,\orcidlink{0000-0002-2997-6731}}
\address[M. P. Saikia]{Mathematical and Physical Sciences Division, School of Arts and Sciences, Ahmedabad University, Navrangpura, Ahmedabad 380009, India}
\email{manjil.saikia@ahduni.edu.in}

\thanks{Corresponding Author: Manjil P. Saikia (manjil.saikia@ahduni.edu.in).}

\keywords{Combinatorial inequalities, Partitions, Parity of parts, Restricted Partitions.}
\subjclass[2020]{05A17, 05A19, 05A20, 11P81.}

\begin{document}

\begin{abstract}
We study inequalities between the eight classical families of integer partitions whose parts are separated by parity in the sense of Andrews. We establish several new comparisons among the counting functions $p_{yx}^{wz}(n)$ and among their non-unitary analogues $q_{yx}^{wz}(n)$, in which parts equal to $1$ are forbidden. All of our proofs are combinatorial. In particular, we show that $q_{od}^{eu}(n)\le q_{eu}^{od}(n)$ for every $n\ge1$, and we determine all cases of equality in both of the conventions used in the literature for partitions whose parts all have the same parity. We also study a parity bias between $q_{od}^{eu}(n)$ and $q_{ed}^{ou}(n)$ in the convention of Ballantine and Welch: in that convention we prove combinatorially that $q_{ed}^{ou}(n)>q_{od}^{eu}(n)$ for all odd $n>7$, and we conjecture that the reverse inequality holds for all even $n>6$.
\end{abstract}

\maketitle

\section{Introduction and statement of results}\label{sec:intro}

A partition of $n$ is a non-increasing sequence $\lambda=(\lambda_1, \lambda_2, \ldots, \lambda_k)$ such that
\[\lambda_1\geq \lambda_2\geq \cdots \geq \lambda_k \quad \text{and}\quad \sum_{i=1}^k\lambda_i = n.\]
For instance the $7$ partitions of $5$ are
\[
5, \quad 4+1, \quad 3+2, \quad 3+1+1, \quad 2+2+1, \quad 2+1+1+1, \quad 1+1+1+1+1.
\]
For a general survey of the theory of partitions we refer the reader to the books of Andrews \cite{gea1} and Johnson \cite{Johnson}, and for the bijective techniques used throughout this paper to the survey of Pak \cite{pak}.

Andrews \cite{andrews-1, andrews} studied partitions in which parts of a given parity are all smaller than those of the other parity, and proved several interesting results, which have since been studied by many other authors (see, for instance, \cite{bringmann, Fu, bringmann2024asymptotic, ballantine2024combinatorial, FX, KS}, and \cite{BCC} for an extension to overpartitions). We denote by $P_{yx}^{wz}(n)$ the classes of partitions of $n$ introduced by Andrews. The symbols $wz$ and $yx$ are each formed with first letter either $e$ or $o$ (denoting even or odd parts) and second letter either $u$ or $d$ (denoting unrestricted or distinct parts). The parts represented by the subscript are assumed to lie below the parts represented by the superscript. This gives rise to eight different families of partitions, namely
\[
P_{eu}^{ou}(n),\quad P_{eu}^{od}(n),\quad P^{eu}_{ou}(n),\quad P_{ou}^{ed}(n),\quad P_{ed}^{ou}(n),\quad P_{ed}^{od}(n),\quad P_{od}^{eu}(n),\quad P_{od}^{ed}(n).
\]
The corresponding cardinalities are denoted by $p_{yx}^{wz}(n)$. Collectively we call all such partitions \emph{partitions with parts separated by parity}.

Two conventions occur in the literature for partitions all of whose parts have a single parity.
\begin{itemize}
\item \textbf{Convention (BW)} (Ballantine--Welch \cite{ballantine2024combinatorial}): partitions all of whose parts are of type $wz$, or all of type $yx$, are included in $P_{yx}^{wz}(n)$.
\item \textbf{Convention (A)} (stated explicitly by Bringmann--Jennings-Shaffer \cite{bringmann}, whose generating functions for all eight families follow this convention, and used by Mahanta--Saikia--Sarma \cite{mahanta2023biases} for the two families with unrestricted parts): every partition in $P_{yx}^{wz}(n)$ must contain at least one part of the lower (subscript) class, \emph{except} in the two families with subscript $eu$, where the even class may be empty (Andrews regards $0$ as the largest even part). Partitions consisting only of parts of the upper class are therefore \emph{excluded} from the six families with subscript $ou$, $od$, $ed$. Thus, for those six families,
\[ p_{yx}^{wz,\,\mathrm{(BW)}}(n) \;=\; p_{yx}^{wz,\,\mathrm{(A)}}(n) \;+\; \#\{\text{partitions of } n \text{ into parts of type } wz\}. \]
\end{itemize}
Each statement below is labelled with the convention to which it applies. For proved inequalities we state explicitly whether the proof works in both conventions; conjectural statements are labelled with the convention in which the numerical evidence supports them.

Andrews \cite{andrews-1} also studied the following subset of $P_{eu}^{ou}(n)$:
\begin{align*}
\overline{P_{eu}^{ou}}(n):=&~\text{the set of partitions of $n$ in $P_{eu}^{ou}(n)$ in which only the largest even part}\\&~\text{appears an odd number of times}.
\end{align*}
Here we allow partitions with no even parts, where we assume that $0$ is the largest even part and that it appears once. We denote the cardinality of this set by $\overline{p_{eu}^{ou}}(n)$, and we note that $\overline{p_{eu}^{ou}}(n)=0$ for odd $n$, since every part other than the largest even part occurs with even multiplicity. Fu and Tang \cite{Fu} considered the analogous set $\overline{P}^{eu}_{ou}(n)$, consisting of those partitions counted by $p^{eu}_{ou}(n)$ in which both even and odd parts occur, the largest even part and the largest odd part each appear with odd multiplicity, and all other parts appear with even multiplicity. We denote its cardinality by $\overline{p}_{ou}^{eu}(n)$.

Bringmann, Craig and Nazaroglu \cite{bringmann2024asymptotic} determined the asymptotic main terms of all eight counting functions. They start from the generating functions of Andrews \cite{andrews} and Bringmann--Jennings-Shaffer \cite{bringmann}, and add to them the contributions of the partitions into parts of a single parity that are omitted there for the six families with subscript $ou$, $od$ and $ed$ (see \cite[footnote~1]{bringmann2024asymptotic}); thus their results are stated in Convention (BW). Their main terms \cite[Theorem~1.1]{bringmann2024asymptotic} imply that in Convention (BW), for $n \gg 1$,
\begin{equation}\label{eq:chain}
	p_{ed}^{od}(n) < p_{od}^{ed}(n) < p_{od}^{eu}(n) < p_{eu}^{od}(n) < p_{ed}^{ou}(n) < p_{eu}^{ou}(n) < p_{ou}^{ed}(n) < p_{ou}^{eu}(n)
\end{equation}
(see also \cite[Section~1]{FX}).
Several of these inequalities were subsequently given combinatorial proofs by Ballantine and Welch \cite{ballantine2024combinatorial}.%; compare also the elementary methods for partition inequalities developed by Berkovich and Uncu \cite{berkovich-uncu}.

\begin{theorem}[Ballantine--Welch \cite{ballantine2024combinatorial}]\label{thm:bw}
    In Convention \textup{(BW)} we have
    \begin{align}
        p_{od}^{ed}(n)&>p_{ed}^{od}(n), & &n\geq 11; \label{A}\\
        p_{ou}^{eu}(n)&>p_{eu}^{ou}(n), & &n\geq 3;\label{bd}\\
        p_{ed}^{ou}(n)&>p_{od}^{eu}(n), & &n\geq 5;\label{e12}\\
        p_{ou}^{ed}(n)&>p_{eu}^{od}(n), & &n\geq 2;\label{e13}\\
        p_{eu}^{od}(n)&>p_{od}^{ed}(n), & &n\geq 8;\label{e14}\\
        \overline{p}_{ou}^{eu}(2n+1)&>\overline{p_{eu}^{ou}}(2n), & &n\geq 3.\label{e15}
    \end{align}
\end{theorem}

\noindent Mahanta, Saikia and Sarma \cite{mahanta2023biases} have also given an elementary proof of \eqref{bd}. Ballantine and Welch \cite{ballantine2024combinatorial} left open combinatorial proofs of the remaining three non-trivial inequalities in the chain of Bringmann, Craig and Nazaroglu, namely
\begin{equation}\label{eq:open}
p_{od}^{eu}(n)<p_{eu}^{od}(n),\qquad p_{eu}^{od}(n)<p_{ed}^{ou}(n),\qquad p_{eu}^{ou}(n)<p_{ou}^{ed}(n).
\end{equation}
The first of these was settled by Fan and Xia \cite{FX}, who gave a combinatorial proof for $n\ge 373$ and checked the range $50\le n\le 400$ by computer.

\begin{theorem}[Fan--Xia \cite{FX}]\label{thm:fx}
In Convention \textup{(BW)}, for all $n\geq 50$ we have
\[p_{eu}^{od}(n)>p_{od}^{eu}(n).\]
\end{theorem}

\noindent The other two inequalities in \eqref{eq:open} still lack combinatorial proofs (see \cite[Section~4]{FX} and \cite[Section~9.3]{KS}).

For our second family of results we consider the \emph{non-unitary} versions of these partitions, that is, partitions in which parts equal to $1$ are not allowed. We denote by $Q_{yx}^{wz}(n)$ the subset of $P_{yx}^{wz}(n)$ consisting of the non-unitary partitions, by $q_{yx}^{wz}(n)$ the corresponding cardinalities, and we retain the two conventions mentioned above. Mahanta, Saikia and Sarma \cite{mahanta2023biases} gave an elementary proof of the following inequality.

\begin{theorem}[Mahanta--Saikia--Sarma \cite{mahanta2023biases}]\label{Qeu}
	In Convention \textup{(A)}, for all $n>3$ we have
	\[q_{ou}^{eu}(n)< q_{eu}^{ou}(n).\]
\end{theorem}

To illustrate the kind of question we consider, take the non-unitary partitions of $9$ with distinct odd parts and unrestricted even parts, in Convention (BW). If the even parts are required to lie above the odd parts, we obtain
\[
Q_{od}^{eu}(9)=\{(9),\ (6,3)\},
\]
while if they are required to lie below the odd parts, we obtain
\[
Q_{eu}^{od}(9)=\{(9),\ (7,2),\ (5,4),\ (5,2,2),\ (3,2,2,2)\}.
\]
The first family embeds into the second: keep $(9)$, and move one unit from the part $3$ to the part $6$ in $(6,3)$, which gives $(7,2)$. Theorem \ref{thm:q4} below shows that placing the odd parts above the even parts never gives fewer partitions, at every size. Its proof is based on the same idea; when several even parts are equal, different amounts have to be moved from different parts, which leads to the staircase transfer of Lemma \ref{lem:transfer}.

The purpose of this paper is to prove the following six inequalities, all by explicit combinatorial injections. Theorems \ref{thm:bar} and \ref{thm:q16} concern the classical (unitary) families, and Theorems \ref{thm:q1}--\ref{thm:q5} concern the non-unitary families. Except for Theorem \ref{thm:q3}, which we deduce from Theorem \ref{thm:q4}, in each case we construct an injection from the smaller family into the larger, together with an explicit element of the larger family that is not in the image.

\begin{theorem}\label{thm:bar}
For all $n>2$, in both conventions, we have
\begin{equation}\label{proved4}
    p_{ou}^{eu}(n)>\overline{p_{eu}^{ou}}(n).
\end{equation}
\end{theorem}

\begin{theorem}\label{thm:q1}
For all $n>3$, in both conventions, we have
\begin{equation}\label{q-1}
q_{eu}^{ou}(n)>q_{od}^{ed}(n).
\end{equation}
\end{theorem}

\begin{theorem}\label{thm:q3}
For all $n>4$, in both conventions, we have
\begin{equation}\label{q-3}
q_{eu}^{ou}(n)>q_{od}^{eu}(n).
\end{equation}
\end{theorem}

\begin{theorem}\label{thm:q4}
In both conventions, for all $n\ge1$ we have
\begin{equation}\label{q-4}
q_{od}^{eu}(n)\le q_{eu}^{od}(n).
\end{equation}
In Convention \textup{(BW)} equality holds exactly for $n\in\{1,2,3,4,6,8\}$, and in Convention \textup{(A)} equality holds exactly for $n\in\{1,3\}$. In particular, $q_{eu}^{od}(n)>q_{od}^{eu}(n)$ for all $n>3$ in Convention \textup{(A)}, and for all $n>8$ in Convention \textup{(BW)}.
\end{theorem}

\begin{theorem}\label{thm:q5}
For all $n>14$, in both conventions, we have
\begin{equation}\label{q-5}
q_{ou}^{ed}(n)>q_{ed}^{od}(n).
\end{equation}
\end{theorem}

\begin{theorem}\label{thm:q16}
For all $n>2$, in both conventions, we have
\begin{equation}\label{q-16}
     p_{ou}^{ed}(n)>p_{ed}^{ou}(n).
\end{equation}
\end{theorem}

\begin{remark}\label{rem:reductions}
The inequalities of Theorems \ref{thm:q1}, \ref{thm:q3} and \ref{thm:q4} are not independent. Indeed, distinct odd upper parts form a subset of unrestricted odd upper parts, so that $q_{eu}^{od}(n)\le q_{eu}^{ou}(n)$, and distinct even upper parts form a subset of unrestricted even upper parts, so that $q_{od}^{ed}(n)\le q_{od}^{eu}(n)$; both inequalities hold in either convention. Hence
\[
q_{od}^{ed}(n)\le q_{od}^{eu}(n)\le q_{eu}^{od}(n)\le q_{eu}^{ou}(n)
\]
for all $n$ in both conventions, and wherever the middle inequality is strict, Theorem \ref{thm:q4} implies Theorems \ref{thm:q3} and \ref{thm:q1}. By Theorem \ref{thm:q4}, this covers all $n>3$ in Convention (A), and $n=5$, $n=7$ and all $n\ge9$ in Convention (BW). The remaining cases in Convention (BW), namely $n\in\{4,6,8\}$ for Theorem \ref{thm:q1} and $n\in\{6,8\}$ for Theorem \ref{thm:q3}, can be read off from Table \ref{tab:smallq}. We prove Theorem \ref{thm:q3} in this way in Section \ref{sec:q4}. Theorem \ref{thm:q1} could be proved in the same way, but we give a direct proof in Section \ref{sec:q1}, since it also yields the intermediate comparison of Proposition \ref{prop:intermediate}; this involves the family $Q_{ed}^{od}(n)$ and does not follow from the chain above.
\end{remark}

\begin{table}[htbp]
\centering
\begin{tabular}{|l|c|c|c|c|c|c|c|c|}
\hline
$n$ & $1$ & $2$ & $3$ & $4$ & $5$ & $6$ & $7$ & $8$ \\ \hline
$q_{eu}^{ou}(n)$ & $0$ & $1$ & $1$ & $2$ & $2$ & $4$ & $3$ & $7$ \\
$q_{eu}^{od}(n)$ & $0$ & $1$ & $1$ & $2$ & $2$ & $3$ & $3$ & $6$ \\ \hline
$q_{od}^{eu}(n)$, (BW) & $0$ & $1$ & $1$ & $2$ & $1$ & $3$ & $2$ & $6$ \\
$q_{od}^{ed}(n)$, (BW) & $0$ & $1$ & $1$ & $1$ & $1$ & $2$ & $2$ & $3$ \\ \hline
$q_{od}^{eu}(n)$, (A) & $0$ & $0$ & $1$ & $0$ & $1$ & $0$ & $2$ & $1$ \\
$q_{od}^{ed}(n)$, (A) & $0$ & $0$ & $1$ & $0$ & $1$ & $0$ & $2$ & $1$ \\ \hline
\end{tabular}
\caption{Small values for Theorems \ref{thm:q1}, \ref{thm:q3} and \ref{thm:q4}. The two families with subscript $eu$ have the same counts in both conventions.}
\label{tab:smallq}
\end{table}

\begin{remark}\label{rem:fx}
Theorem \ref{thm:fx} has the following consequence for Convention (A). Since the subscript of $P_{eu}^{od}$ is $eu$, we have $p_{eu}^{od,\mathrm{(A)}}(n)=p_{eu}^{od,\mathrm{(BW)}}(n)$, whereas $p_{od}^{eu,\mathrm{(A)}}(n)\le p_{od}^{eu,\mathrm{(BW)}}(n)$, with equality for odd $n$. Hence $p_{eu}^{od}(n)>p_{od}^{eu}(n)$ holds in both conventions for all $n\ge 50$. A direct computation for $n<50$ then shows that, in Convention (A), $p_{eu}^{od}(n)>p_{od}^{eu}(n)$ for \emph{every} even $n\ge 2$. On the other hand, the inequality fails in both conventions for every odd $n<50$, and in Convention (BW) it fails for every even $n\le 16$.
\end{remark}

If instead the parts in the lower block are required to be distinct, the comparison depends on the parity of $n$. In Convention (BW), the family $Q_{od}^{eu}(9)$ above has two elements, while
\[
Q_{ed}^{ou}(9)=\{(9),\ (7,2),\ (5,4),\ (3,3,3)\}
\]
has four: the odd parts may now repeat, but the even parts may not. For $n=8$ the comparison goes the other way, since $q_{od}^{eu}(8)=6$ and $q_{ed}^{ou}(8)=4$. We prove one half of this parity bias.

\begin{theorem}\label{thm:parityodd}
In Convention \textup{(BW)}, for all odd $n>7$ we have
\begin{equation}\label{q-odd}
q_{ed}^{ou}(n)>q_{od}^{eu}(n).
\end{equation}
\end{theorem}

\noindent Numerical data suggest that the inequality is reversed for even $n>6$; see Conjecture \ref{conj:parity}. We also prove a quantitative strengthening of Theorem \ref{thm:q16} in Convention (A).

\begin{theorem}\label{thm:ratio2}
In Convention \textup{(A)}, for all $n\ge3$ we have
\begin{equation}\label{q-ratio}
p_{ou}^{ed}(n)>2\,p_{ed}^{ou}(n).
\end{equation}
\end{theorem}

\noindent In Section \ref{sec:conclusion} we state Conjecture \ref{conj:parity} and some further conjectures suggested by numerical data.

The rest of the paper is organised as follows. Section \ref{sec:transfer} contains a lemma used in several of the proofs, Section \ref{sec:pf} contains the proof of Theorem \ref{thm:bar}, Section \ref{sec:q1} contains the proof of Theorem \ref{thm:q1}, Section \ref{sec:q4} contains the proofs of Theorems \ref{thm:q4} and \ref{thm:q3}, Section \ref{sec:q5} contains the proof of Theorem \ref{thm:q5}, Section \ref{sec:q16} contains the proof of Theorem \ref{thm:q16}, Section \ref{sec:parity} contains the proof of Theorem \ref{thm:parityodd}, and Section \ref{sec:ratio} contains the proof of Theorem \ref{thm:ratio2}. We close with some conjectures and open problems in Section \ref{sec:conclusion}.

\subsection*{Human disclosure} This project started almost half a decade ago, when the last author was working on \cite{mahanta2023biases}. He made several observations some of which are now theorems in this paper. In 2024 he started a different collaboration with the first author when he mentioned these observations and shared some conjectures. After some back and forth, the two of them made some partial progress. When the work stalled, the last author asked the third author to join them, sometime in 2025. After some further progress, the third author recruited the second author into the project who made further progress until the first half of 2026. All of this was done without the use of AI. After this, the project stalled for several months, until November 2026, when we used AI, which we now disclose below.

\subsection*{AI disclosure}
Parts of this work were carried out with the assistance of Claude (Claude Opus 5.5, developed by Anthropic). Specifically, AI was used for the following.
\begin{itemize}
\item \emph{Checking the manuscript.} Claude computed all of the counting functions considered here (up to $n=400$), checked every stated inequality, bound and table of small values against these data, and tested by computer the injections in the proofs of Theorems \ref{thm:q1}, \ref{thm:q3}, \ref{thm:q4}, \ref{thm:q5} and \ref{thm:q16} (for $n$ up to between $36$ and $125$, depending on the theorem). This led to a number of corrections. These include the convention-dependent statements in Theorem \ref{thm:bar} and in the conjectures of Section \ref{sec:conclusion}, and the observation that an inequality we had stated as a conjecture follows from the work of Fan and Xia \cite{FX} (Remark \ref{rem:fx}).
\item \emph{Proofs.} Claude proposed the proof of Theorem \ref{thm:bar} in both conventions, Theorem \ref{thm:parityodd} and its proof, and Theorem \ref{thm:ratio2} together with Lemma \ref{lem:Phi} and its proof. In each case the arguments were also checked by computer for small $n$.
\item \emph{Referee-style reports.} Two referee-style reports on an earlier version of this paper were generated with ChatGPT (models GPT-6.2 Astra and ChatGPT 6.2 Sol, developed by OpenAI). The resulting revisions were made with Claude's assistance. They include the extension of Theorem \ref{thm:q4} to Convention (BW) together with its proof (suggested in the first report), the correction identity \eqref{eq:correction} in Section \ref{sec:conclusion}, the treatment of several boundary cases in the proofs of Theorem \ref{thm:parityodd} and Lemma \ref{lem:Phi}, Proposition \ref{prop:limits} and its proof (suggested in the second report), and the proof of Theorem \ref{thm:q5} in Section \ref{sec:q5}.
\item \emph{Editing.} Claude assisted with the exposition, drafted the remarks and conjectures in Section \ref{sec:conclusion}, and restructured the proofs in Sections \ref{sec:q1}, \ref{sec:q4} and \ref{sec:q16} around Lemma \ref{lem:transfer}, which isolates the verification common to the injections constructed there.
\item \emph{Drawing.} All diagrams were generated by Claude.
\item \emph{Final referee pass.} Claude Fable 5.2 was used for a final referee pass. No major corrections were suggested in this pass\footnote{Some suggestions by Claude and ChatGPT were willfully ignored to give precedence to clarity over brevity.}.
\end{itemize}
The authors have reviewed all of this material and take full responsibility for the contents of this paper.

\section{A transfer lemma}\label{sec:transfer}

Several of our injections convert a mixed partition by exchanging the parities of some of its largest parts and some of its smallest parts. Consider, for instance, $(6,6,5,3)\in Q_{od}^{eu}(20)$. Adding $1$ to both parts $6$ and subtracting $1$ from both odd parts gives $(7,7,4,2)$, whose odd parts are no longer distinct. Adding $3$ and $1$ to the two parts $6$ instead produces the distinct odd parts $9$ and $7$, and subtracting $3$ and $1$ from the odd parts $5$ and $3$ pays for these increases; the result $(9,7,2,2)$ lies in $Q_{eu}^{od}(20)$, where the even parts may repeat. The following lemma records the general form of this operation, together with the simpler variant in which every shift equals $1$.

\begin{lemma}\label{lem:transfer}
Let $U=(u_1,\ldots,u_a)$ with $u_1\ge\cdots\ge u_a$ be positive integers of one parity, and let $D=(d_1,\ldots,d_b)$ with $d_1>\cdots>d_b\ge2$ be integers of the other parity, where $a,b\ge1$ and $u_a>d_1$. Let $r=\min(a,b)$, and let $(s_1,\ldots,s_r)$ be one of the following two sequences:
\begin{itemize}
\item[(S)] the \emph{staircase} $s_i=2(r-i)+1$ for $1\le i\le r$;
\item[(O)] the \emph{constant} sequence $s_i=1$ for $1\le i\le r$.
\end{itemize}
Let $T(U,D)$ be the partition whose parts are the \emph{upper parts}
\[
u_1+s_1,\ \ldots,\ u_r+s_r,\ d_1,\ \ldots,\ d_{b-r}
\]
together with the \emph{lower parts}
\[
u_{r+1},\ \ldots,\ u_a,\ d_{b-r+1}-s_1,\ \ldots,\ d_b-s_r.
\]
Then the following hold.
\begin{enumerate}
\item[(a)] $|T(U,D)|=|U|+|D|$. The $b$ upper parts have the parity of $D$, and the $a$ lower parts have the parity of $U$. Every upper part exceeds every lower part, and the smallest part of $T(U,D)$ is $d_b-1$.
\item[(b)] In the order displayed, the lower parts are weakly decreasing, and the upper parts are weakly decreasing. In case \textup{(S)} the upper parts are distinct. In case \textup{(O)} the parts $d_{b-r+i}-1$ are distinct, and if the parts of $U$ are distinct, then all parts of $T(U,D)$ are distinct.
\item[(c)] $T(U,D)$ (together with the choice of \textup{(S)} or \textup{(O)}) determines $(U,D)$: the numbers $a$ and $b$ are the numbers of parts of $T(U,D)$ of the parity of $U$ and of $D$, respectively; with $r=\min(a,b)$, the $r$ largest parts of the parity of $D$ are $u_1+s_1\ge\cdots\ge u_r+s_r$, the $r$ smallest parts of the parity of $U$ are $d_{b-r+1}-s_1\ge\cdots\ge d_b-s_r$, and the remaining parts are $u_{r+1},\ldots,u_a$ and $d_1,\ldots,d_{b-r}$.
\end{enumerate}
\end{lemma}

\begin{proof}
(a) The shifts are added to $r$ parts and subtracted from $r$ parts, so the size is unchanged, and since each $s_i$ is odd, each converted part changes parity. Note that $a=r$ or $b=r$. Since $u_r+s_r\ge u_a+1>d_1$, the smallest upper part is $u_r+1$ if $b=r$, and $d_{b-r}$ if $b>r$. The largest lower part is $u_{r+1}$ if $a>r$ (and then $b=r$, so $u_{r+1}\ge u_a>d_1>d_1-s_1$), and $d_{b-r+1}-s_1$ otherwise. Comparing these in the three cases $a>r$, $b>r$ and $a=b=r$ shows that every upper part exceeds every lower part: $u_r+1>u_{r+1}$, $d_{b-r}>d_{b-r+1}-s_1$ and $u_r+1>d_1-s_1$, respectively. The smallest part is the last lower part $d_b-s_r=d_b-1$.

(b) For $1\le i<r$ we have $(u_i+s_i)-(u_{i+1}+s_{i+1})=(u_i-u_{i+1})+(s_i-s_{i+1})$ and $(d_{b-r+i}-s_i)-(d_{b-r+i+1}-s_{i+1})=(d_{b-r+i}-d_{b-r+i+1})-(s_i-s_{i+1})$, where $s_i-s_{i+1}$ equals $2$ in case (S) and $0$ in case (O), and $d_{b-r+i}-d_{b-r+i+1}\ge2$. Together with $u_r+s_r>d_1$ and, when $a>r$, $u_a>d_1-s_1$, this gives all the stated monotonicity and distinctness properties.

(c) By (a), the upper parts are exactly the parts of $T(U,D)$ of the parity of $D$, and the lower parts are exactly those of the parity of $U$. By (b), listing them as displayed lists each block in weakly decreasing order, so the converted parts are the $r$ largest parts of the first block and the $r$ smallest parts of the second block. Subtracting, respectively adding, $s_1,\ldots,s_r$ recovers $U$ and $D$.
\end{proof}

For example, with $U=(6,6,6)$ and $D=(5,3)$ we have $r=2$; the staircase $(3,1)$ gives $T(U,D)=(9,7,6,2,2)$, while the constant sequence gives $T(U,D)=(7,7,6,4,2)$. In our applications $d_b\ge3$ for the non-unitary families, so that $T(U,D)$ has no part equal to $1$, and $d_b\ge2$ for the unitary families.

Throughout Sections \ref{sec:q1}--\ref{sec:ratio}, $k$ denotes the number of odd parts and $\ell$ the number of even parts of the partition under consideration.

\section{Proof of Theorem \ref{thm:bar}}\label{sec:pf}

\begin{proof}[Proof of Theorem \ref{thm:bar}]
Note first that the family $\overline{P_{eu}^{ou}}(n)$ is the same in both conventions, and that every partition of $n$ into odd parts lies in $P_{ou}^{eu}(n)$ in both conventions (its lower class is non-empty). We construct an injection $\varphi\colon\overline{P_{eu}^{ou}}(n)\to P_{ou}^{eu,\mathrm{(A)}}(n)$; since $P_{ou}^{eu,\mathrm{(A)}}(n)\subseteq P_{ou}^{eu,\mathrm{(BW)}}(n)$, this proves \eqref{proved4} in both conventions.

If $n$ is odd, then $\overline{p_{eu}^{ou}}(n)=0$, while $(n)\in P_{ou}^{eu}(n)$, so \eqref{proved4} holds. Let now $n\ge4$ be even and $\lambda\in\overline{P_{eu}^{ou}}(n)$. Let $e$ be the largest even part of $\lambda$, with $e=0$ if $\lambda$ has no even part. Since every part of $\lambda$ other than one copy of $e$ occurs with even multiplicity, there is a unique partition $\mu=(\mu_1,\ldots,\mu_r)$ such that $\lambda=(e)\cup\mu\cup\mu$ (where $(e)$ is omitted if $e=0$); the odd parts of $\mu$ exceed $e$ and its even parts are at most $e$. Define
\[
\varphi(\lambda):=\begin{cases}
\lambda, & \text{if } e=0;\\
(n-1,1), & \text{if } \lambda=(n);\\
(2\mu_1,\ldots,2\mu_r,\underbrace{1,\ldots,1}_{e}), & \text{otherwise}.
\end{cases}
\]
In the first case $\varphi(\lambda)$ is a partition into odd parts. In the second case $n-1\ge3$ is odd, so again $\varphi(\lambda)$ consists of odd parts only. In the third case $r\ge1$ and $e\ge2$, so $\varphi(\lambda)$ contains at least one part equal to $1$ and at least one even part, and its odd parts (all equal to $1$) lie below its even parts; moreover $|\varphi(\lambda)|=2|\mu|+e=n$. Hence $\varphi(\lambda)\in P_{ou}^{eu,\mathrm{(A)}}(n)$ in all cases.

The three cases have disjoint images: images in the first case consist of odd parts, each occurring with even multiplicity; the image $(n-1,1)$ of the second case consists of odd parts and contains the part $n-1$ with multiplicity one; images in the third case contain an even part. In the third case $\lambda$ is recovered from $\varphi(\lambda)$, since $e$ is the number of parts equal to $1$ and $\mu$ is obtained by halving the even parts. Thus $\varphi$ is injective.

Finally, $\omega_n=(2,1^{n-2})\in P_{ou}^{eu,\mathrm{(A)}}(n)$ is not in the image of $\varphi$. Indeed, it contains an even part, so it could only arise in the third case with $e=n-2$ and $\mu=(1)$; but then $\lambda=(n-2,1,1)$ would have the even part $n-2\ge 2$ above the odd part $1$, contradicting $\lambda\in P_{eu}^{ou}(n)$. This proves \eqref{proved4} for all $n\ge3$.

The bound $n>2$ is sharp: $p_{ou}^{eu}(2)=2=\overline{p_{eu}^{ou}}(2)$ in Convention (BW), both sets being $\{(2),(1,1)\}$, and $p_{ou}^{eu}(2)=1<2=\overline{p_{eu}^{ou}}(2)$ in Convention (A).
\end{proof}

\begin{remark}
In Convention (BW), Theorem \ref{thm:bar} can also be deduced from \eqref{bd}: since $\overline{P_{eu}^{ou}}(n)\subseteq P_{eu}^{ou}(n)$, we have $p_{ou}^{eu}(n)>p_{eu}^{ou}(n)\ge\overline{p_{eu}^{ou}}(n)$ for $n\ge3$. This argument does not apply in Convention (A), where \eqref{bd} fails for $n=4$ and $n=6$.
\end{remark}

\section{Proof of Theorem \ref{thm:q1}}\label{sec:q1}

The proof of Theorem \ref{thm:q1} passes through the family $Q_{ed}^{od}(n)$, in which both blocks consist of distinct parts, and gives the following intermediate comparison.

\begin{proposition}\label{prop:intermediate}
In Convention \textup{(BW)}, for all $n\ge9$ we have
\[
q_{od}^{ed}(n)<q_{ed}^{od}(n)<q_{eu}^{ou}(n).
\]
\end{proposition}

\begin{proof}[Proof of Proposition \ref{prop:intermediate} and Theorem \ref{thm:q1}]
The proof is carried out in Convention (BW). Since $q_{yx}^{wz,\mathrm{(A)}}\le q_{yx}^{wz,\mathrm{(BW)}}$, with equality for the family $Q_{eu}^{ou}$, the Convention (A) version of \eqref{q-1} follows as well.

We construct an injection proving $q_{od}^{ed}(n)<q_{ed}^{od}(n)<q_{eu}^{ou}(n)$ for $n\ge9$; the cases $4\le n\le 8$ are checked directly at the end.

If $\lambda\in Q^{ed}_{od}(n)$ consists of either distinct odd parts only or distinct even parts only, then clearly $\lambda\in Q^{od}_{ed}(n)$, and we let $\psi(\lambda)=\lambda$. For the remaining (mixed) partitions $\lambda=(e_1,\ldots,e_\ell,o_1,\ldots,o_k)\in Q^{ed}_{od}(n)$, with distinct even parts $e_1>\cdots>e_\ell$, distinct odd parts $e_\ell>o_1>\cdots>o_k\ge 3$, and $k,\ell\ge1$, let $\psi(\lambda):=T(U,D)$ with $U=(e_1,\ldots,e_\ell)$, $D=(o_1,\ldots,o_k)$ and the constant sequence (O) in Lemma \ref{lem:transfer}. Explicitly,
\[
\psi(\lambda)=
\begin{cases}
    (e_1+1,\ldots,e_\ell+1,\,o_1,\ldots,o_{k-\ell},\,o_{k-\ell+1}-1,\ldots,o_k-1), & \text{ if } \ell\le k;\\
    (e_1+1,\ldots,e_k+1,\,e_{k+1},\ldots,e_\ell,\,o_1-1,\ldots,o_k-1), & \text{ if } \ell>k.
\end{cases}
\]
Since the parts of $U$ are distinct, Lemma \ref{lem:transfer} shows that $\psi(\lambda)$ is a partition of $n$ into $k$ odd parts lying above $\ell$ even parts, with all parts distinct and the smallest part equal to $o_k-1\ge2$. Hence $\psi(\lambda)\in Q^{od}_{ed}(n)$. Pure partitions are fixed by $\psi$ and have parts of one parity only, whereas the images of mixed partitions contain parts of both parities, so by Lemma \ref{lem:transfer}(c) the map $\psi$ is injective. Hence $q_{od}^{ed}(n)\le q_{ed}^{od}(n)$ for all $n$.

\noindent\emph{Strictness for $n\ge9$.} The following partition lies in $Q^{od}_{ed}(n)$ and is not in the image of $\psi$:
$$
\omega_n=\begin{cases}
\left(\frac{n+1}{2},\frac{n-1}{2}\right), & \text{ if } n\equiv 1\pmod 4;\\[2pt]
\left(\frac{n-1}{2},\frac{n-3}{2},2\right), & \text{ if } n \equiv 3 \pmod 4;\\[2pt]
\left(n-5,3,2\right), & \text{ if } n \text{ is even}.
\end{cases}
$$
Indeed, $\omega_n$ is mixed, so it is not a fixed pure partition, and a preimage $\lambda$ would have as many even parts $\ell$ and odd parts $k$ as $\omega_n$; we apply the inverse described in Lemma \ref{lem:transfer}(c), with $r=1$ in each case. If $n\equiv1\pmod4$, then $\ell=k=1$, and the inverse would give $e_1=\frac{n+1}{2}-1=\frac{n-1}{2}$ and $o_1=\frac{n-1}{2}+1=\frac{n+1}{2}>e_1$, contradicting $e_1>o_1$. If $n\equiv3\pmod4$, then $\ell=2$ and $k=1$, and the inverse would give $e_1=\frac{n-1}{2}-1=\frac{n-3}{2}$, $o_1=2+1=3$ and the retained part $e_2=\frac{n-3}{2}=e_1$, contradicting distinctness. If $n$ is even, then $\ell=1$ and $k=2$, and the inverse would give $e_1=n-6$, the retained part $o_1=3$ and $o_2=2+1=3=o_1$, again contradicting distinctness. For $n\ge 9$ all parts of $\omega_n$ are $\ge2$, they are distinct (when $n\equiv3\pmod4$ we have $n\ge11$, so $\frac{n-3}{2}>2$), and the parity separation holds, so $\omega_n\in Q_{ed}^{od}(n)$. Hence $q_{od}^{ed}(n)<q_{ed}^{od}(n)$ for $n\ge9$.

Next, $Q^{od}_{ed}(n)\subseteq Q^{ou}_{eu}(n)$, since relaxing both distinctness conditions preserves membership. Moreover, for $n\ge9$ the partition
$$\begin{cases}
    (n-4,2,2), & \text{ if } n \text{ is odd};\\
    (n-7,3,2,2), & \text{ if } n \text{ is even},
\end{cases}
$$
lies in $Q^{ou}_{eu}(n)$ but not in $Q^{od}_{ed}(n)$, because the even part $2$ is repeated. Hence
$$q^{ed}_{od}(n)< q^{od}_{ed}(n)< q^{ou}_{eu}(n),\qquad n\ge 9,$$
which is Proposition \ref{prop:intermediate} and gives \eqref{q-1} for $n\ge9$. For the remaining cases one checks directly that
\[ \big(q_{eu}^{ou}(n),\,q_{od}^{ed}(n)\big)_{n=4,\dots,8}=(2,1),\,(2,1),\,(4,2),\,(3,2),\,(7,3), \]
completing the proof of \eqref{q-1} for all $n>3$.
\end{proof}

\section{Proofs of Theorems \ref{thm:q4} and \ref{thm:q3}}\label{sec:q4}

\begin{proof}[Proof of Theorem \ref{thm:q4}]
We first work in Convention (BW), and deduce the statement for Convention (A) at the end. We construct an injection $\psi\colon Q_{od}^{eu}(n)\to Q_{eu}^{od}(n)$.

If $\lambda \in Q_{od}^{eu}(n)$ consists of distinct odd parts only, or of even parts only, then clearly $\lambda \in Q_{eu}^{od}(n)$, and we let $\psi(\lambda)=\lambda$. For the remaining (mixed) partitions
\[\lambda = (\lambda_{e_1},\ldots,\lambda_{e_\ell},\lambda_{o_1},\ldots,\lambda_{o_k})\in Q_{od}^{eu}(n)\]
with unrestricted even parts $\lambda_{e_1} \geq \cdots \geq \lambda_{e_\ell}$, distinct odd parts $\lambda_{o_1} > \cdots > \lambda_{o_k}\geq 3$, and $k,\ell\geq 1$, let $\psi(\lambda):=T(U,D)$ with $U=(\lambda_{e_1},\ldots,\lambda_{e_\ell})$, $D=(\lambda_{o_1},\ldots,\lambda_{o_k})$ and the staircase (S) in Lemma \ref{lem:transfer}. Explicitly, with $r=\min(\ell,k)$,
\begin{multline*}
\psi(\lambda)=\big(\lambda_{e_1}+2r-1,\,\lambda_{e_2}+2r-3,\,\ldots,\,\lambda_{e_r}+1,\ \lambda_{e_{r+1}},\ldots,\lambda_{e_\ell},\\
\lambda_{o_1},\ldots,\lambda_{o_{k-r}},\ \lambda_{o_{k-r+1}}-(2r-1),\,\ldots,\,\lambda_{o_k}-1\big),
\end{multline*}
where at most one of the two groups of unchanged parts is non-empty. By Lemma \ref{lem:transfer}, $\psi(\lambda)$ is a partition of $n$ into $k$ distinct odd parts lying above $\ell$ even parts, with smallest part $\lambda_{o_k}-1\ge2$, so $\psi(\lambda)\in Q_{eu}^{od}(n)$. Pure partitions are fixed by $\psi$ and have parts of only one parity, whereas the images of mixed partitions contain parts of both parities; hence $\psi$ is injective by Lemma \ref{lem:transfer}(c).

\medskip
\noindent\textbf{Example.} Let $n=26$ and $\lambda=(6,6,6,5,3)$, where $\ell=3$, $k=2$ and $r=2$. Then $\psi(\lambda)=(6+3,6+1,6,5-3,3-1)=(9,7,6,2,2)$. To invert, subtract $3$ and $1$ from the two largest odd parts and add $3$ and $1$ to the two smallest even parts; the repeated parts $2$ cause no ambiguity.
\vspace{0.3cm}

\ytableausetup{mathmode, boxsize=1.2em}
\begin{center}
    \begin{ytableau}
        ~ & ~ & ~ & ~ & ~ & ~ \\
        ~ & ~ & ~ & ~ & ~ & ~ \\
        ~ & ~ & ~ & ~ & ~ & ~ \\
        *(gray!30) & *(gray!30) & *(gray!30) & *(gray!30) & *(gray!30) \\
        *(gray!30) & *(gray!30) & *(gray!30)
    \end{ytableau}
    \quad $\longrightarrow$ \quad
    \begin{ytableau}
        ~ & ~ & ~ & ~ & ~ & ~ & *(gray!30) & *(gray!30) & *(gray!30)\\
         ~ & ~ & ~ & ~ & ~ & ~ & *(gray!30) \\
        ~ & ~ & ~ & ~ & ~ & ~ \\
       *(gray!30) & *(gray!30)   \\
        *(gray!30) & *(gray!30)
    \end{ytableau}
\end{center}

\noindent\emph{Strictness for odd $n\ge5$.} The partition
\[
\omega_n=(3,\underbrace{2,2,\ldots,2}_{\frac{n-3}{2}})
\]
lies in $Q_{eu}^{od}(n)$ and is not in the image of $\psi$. Indeed, it is mixed, so it is not a fixed pure partition. It has $\frac{n-3}{2}\ge1$ even parts and one odd part, so a preimage would have $\ell=\frac{n-3}{2}$, $k=1$ and $r=1$, and the inverse of Lemma \ref{lem:transfer}(c) would give $\lambda_{e_1}=3-1=2$ and $\lambda_{o_1}=2+1=3$, contradicting $\lambda_{e_\ell}>\lambda_{o_1}$. (Equivalently, the largest part of every mixed image is an odd part greater than $\lambda_{e_1}\ge4$, since $\lambda_{e_1}>\lambda_{o_1}\ge3$; so it is at least $5$.)

\noindent\emph{Strictness for even $n\ge10$.} Let $t=\frac{n-8}{2}\ge1$. The partition
\[
\omega_n=(5,3,\underbrace{2,2,\ldots,2}_{t})
\]
lies in $Q_{eu}^{od}(n)$ and is not in the image of $\psi$. Indeed, it is mixed, with $t$ even parts and two odd parts, so a preimage would have $\ell=t$ and $k=2$. If $t=1$, then $r=1$: the inverse retains the odd part $3=\lambda_{o_1}$ and recovers $\lambda_{o_2}=2+1=3$ from the even part, contradicting the distinctness of the odd parts. If $t\ge2$, then $r=2$, and the largest odd part of the image is $\lambda_{e_1}+3$; since $\lambda_{e_1}>\lambda_{o_1}>\lambda_{o_2}\ge3$, we have $\lambda_{e_1}\ge6$, so this part is at least $9>5$.

Hence $q_{od}^{eu}(n)<q_{eu}^{od}(n)$ in Convention (BW) for $n=5$, $n=7$ and all $n\ge9$. For the remaining values $n\in\{1,2,3,4,6,8\}$ the two counts coincide in Convention (BW); they are $0,1,1,2,3,6$, respectively (see Table \ref{tab:smallq}). For instance, for $n=6$ both families consist of $(6)$, $(4,2)$ and $(2,2,2)$. This proves the statement for Convention (BW).

\noindent\emph{Convention \textup{(A)}.} Since the subscript of $Q_{eu}^{od}$ is $eu$, we have $Q_{eu}^{od,\mathrm{(A)}}(n)=Q_{eu}^{od,\mathrm{(BW)}}(n)$, whereas $Q_{od}^{eu,\mathrm{(A)}}(n)$ is obtained from $Q_{od}^{eu,\mathrm{(BW)}}(n)$ by removing the partitions of $n$ into even parts. There are $p(n/2)\ge1$ of these for even $n$, where $p$ denotes the ordinary partition function, and none for odd $n$. Hence, in Convention (A), the inequality $q_{od}^{eu}(n)\le q_{eu}^{od}(n)$ holds for all $n$, it is strict for every even $n\ge2$ and for every odd $n\ge5$, and equality holds for $n\in\{1,3\}$, where both counts are $0$ and $1$, respectively. (In terms of the map $\psi$, which in Convention (A) has no pure even partitions to fix, the pure even partition $(n)$ is not in the image for even $n$.)
\end{proof}

\begin{proof}[Proof of Theorem \ref{thm:q3}]
Relaxing the distinctness condition on the odd parts gives $Q_{eu}^{od}(n)\subseteq Q_{eu}^{ou}(n)$ in both conventions. Together with Theorem \ref{thm:q4} this gives
\[
q_{od}^{eu}(n)\le q_{eu}^{od}(n)\le q_{eu}^{ou}(n)\qquad\text{for all } n\ge1,
\]
in both conventions, and the first inequality is strict for all $n\ge4$ in Convention (A), and for $n=5$, $n=7$ and all $n\ge9$ in Convention (BW). For the two remaining values in Convention (BW), Table \ref{tab:smallq} gives $q_{od}^{eu}(6)=3<4=q_{eu}^{ou}(6)$ and $q_{od}^{eu}(8)=6<7=q_{eu}^{ou}(8)$. This proves \eqref{q-3} for all $n>4$ in both conventions. The bound $n>4$ is sharp in Convention (BW), since $q_{od}^{eu}(n)=q_{eu}^{ou}(n)$ for $1\le n\le4$ (see Table \ref{tab:smallq}).
\end{proof}

\section{Proof of Theorem \ref{thm:q5}}\label{sec:q5}

To see why the proof of Theorem \ref{thm:q5} requires a different idea from those of Sections \ref{sec:q1} and \ref{sec:q4}, consider $\lambda=(9,8,4,2)\in Q_{ed}^{od}(23)$. Moving one unit from the part $2$ to the part $9$ produces a part equal to $1$, which is not allowed, while moving one unit from the part $4$ instead produces the odd part $3$ lying between the even parts $8$ and $2$. We shall map $\lambda$ to $(10,7,3,3)\in Q_{ou}^{ed}(23)$, which has three odd parts although $\lambda$ has only one. To allow the numbers of parts of each parity to change without losing control of the inverse, we record each partition relative to the smallest partition with the same numbers of odd and even parts.

\begin{proof}[Proof of Theorem \ref{thm:q5}]

\noindent\emph{Minimum partitions and gap coordinates.} Throughout this proof, minimum partitions are written with their parts in \emph{increasing} order, and the coordinates defined below refer to this order. For $k,\ell\ge1$, the smallest partitions with $k$ odd parts and $\ell$ even parts in $Q_{ed}^{od}$ and in $Q_{ou}^{ed}$ are, respectively,
\begin{align*}
A_{k,\ell}&:=(2,4,\ldots,2\ell,\ 2\ell+1,2\ell+3,\ldots,2\ell+2k-1),\\
B_{k,\ell}&:=(\underbrace{3,\ldots,3}_{k},\ 4,6,\ldots,2\ell+2).
\end{align*}
For pure partitions, the smallest partitions with $\ell$ even parts are $(2,4,\ldots,2\ell)$ in both families, while the smallest partitions with $k$ odd parts are $(3,5,\ldots,2k+1)$ in $Q_{ed}^{od}$ and $(3,\ldots,3)$ ($k$ times) in $Q_{ou}^{ed}$.

Let $M=(m_1,\ldots,m_L)$ be one of these minimum partitions, of length $L$. Then every partition in the corresponding family with the same numbers of odd and even parts can be written uniquely as $(m_1+2y_1,\ldots,m_L+2y_L)$ with integers $0\le y_1\le\cdots\le y_L$. Indeed, consecutive parts of the same parity differ by at least $2$ if they are required to be distinct and by at least $0$ otherwise, and at the change of parity they differ by at least $1$; these are exactly the differences between consecutive parts of $M$. Put $y_0=0$ and define the \emph{gap coordinates}
\[
g_j:=y_{L-j+1}-y_{L-j}\ge0,\qquad 1\le j\le L.
\]
Increasing $g_j$ by $1$ increases each of the $j$ largest parts by $2$, so that
\begin{equation}\label{eq:gapweight}
|\lambda|=|M|+2\sum_{j=1}^{L}j\,g_j.
\end{equation}
Conversely, any vector of non-negative integers $(g_1,\ldots,g_L)$ arises in this way: $y_i=g_{L-i+1}+\cdots+g_L$. In particular, $g_1$ measures the excess of the difference between the two largest parts over its minimum value, divided by $2$, and $g_L=y_1$ measures the excess of the smallest part over $m_1$, divided by $2$.

For example, $\lambda=(11,9,8,4)$ has $k=\ell=2$; subtracting $A_{2,2}=(2,4,5,7)$ from $(4,8,9,11)$ and halving gives $y=(1,2,2,2)$, so $g=(0,0,1,1)$.

\medskip
\noindent\emph{The map on mixed partitions.} Let $n\ge20$, and let $\lambda\in Q_{ed}^{od}(n)$ be mixed, with gap coordinates $g=(g_1,\ldots,g_{k+\ell})$ relative to $A_{k,\ell}$. We define $\psi(\lambda)\in Q_{ou}^{ed}(n)$ by prescribing its numbers $(k',\ell')$ of odd and even parts and its gap coordinates $h$ relative to the corresponding minimum partition $B_{k',\ell'}$ (or $(3,3,3)$ if $\ell'=0$). Unless stated otherwise, $h_j=g_j$. Put
\[
D_{k,\ell}:=\frac{|A_{k,\ell}|-|B_{k,\ell}|}{2}=(k-1)\ell+\frac{k(k-3)}{2}.
\]
Note that $D_{1,\ell}=-1$, $D_{2,\ell}=\ell-1\ge0$ and $D_{k,\ell}\ge0$ for $k\ge3$.
\begin{enumerate}
\item[(i)] If $k\ge2$, or $k=1$ and $g_1>0$: let $(k',\ell')=(k,\ell)$ and $h_1=g_1+D_{k,\ell}$.
\item[(ii)] If $k=1$, $g_1=0$ and $\ell\ge3$: let $(k',\ell')=(3,\ell-2)$, $h_1=2\ell-5$ and $h_2=g_2+1$.
\item[(iii)] If $k=1$, $g_1=0$ and $\ell=2$: let $(k',\ell')=(3,0)$ and $h_1=1$.
\item[(iv)] If $k=1$, $g_1=0$ and $\ell=1$: let $(k',\ell')=(3,0)$ and $h=(0,\,g_2-1,\,0)$.
\end{enumerate}
In each case all coordinates of $h$ are non-negative: in (i) this is clear if $k\ge2$, and $h_1=g_1-1\ge0$ if $k=1$; in (iv) we have $n=|(2,3)|+4g_2=5+4g_2$, so $g_2\ge4$. In cases (i)--(iii) the output has the same length $k+\ell$ as $\lambda$, and in case (iv) the length increases from $2$ to $3$. By \eqref{eq:gapweight} each rule preserves the size: in (i) the change $2D_{k,\ell}$ in $h_1$ compensates the difference $|A_{k,\ell}|-|B_{k,\ell}|$; in (ii) the minimum partition decreases by $|A_{1,\ell}|-|B_{3,\ell-2}|=4\ell-6$, which is restored by $2(2\ell-5)+2\cdot2=4\ell-6$; in (iii) the minimum decreases from $2+4+5=11$ to $3+3+3=9$, which is restored by $2h_1=2$; and in (iv) the minimum increases from $2+3=5$ to $9$, which is compensated by the decrease of $g_2$ by $1$ at weight $2\cdot2$. Hence $\psi(\lambda)\in Q_{ou}^{ed}(n)$ in both conventions, since it contains an odd part.

The images of the mixed partitions are distinguished as follows. A mixed target with $\ell'\ge1$ even parts arises from rule (i), or from rule (ii) when $k'=3$. In the latter case $\ell'=\ell-2$ and $h_1=2\ell-5=2\ell'-1$, while rule (i) with $(k,\ell)=(3,\ell')$ gives $h_1\ge D_{3,\ell'}=2\ell'$. A pure odd target arises from a mixed partition only through rules (iii) and (iv); these targets have exactly three parts, and they have $h_1=1$ and $h_1=0$, respectively. Once the rule is known, the numbers $(k,\ell)$ and the coordinates $g$ are recovered by undoing the stated changes; for instance, in rule (ii) we have $\ell=\ell'+2$, $g_1=0$ and $g_2=h_2-1$. Thus $\psi$ is injective on mixed partitions.

For example, $\lambda=(11,9,8,4)$ above falls under rule (i) with $D_{2,2}=1$, so $h=(1,0,1,1)$ relative to $B_{2,2}=(3,3,4,6)$, and $\psi(\lambda)=(12,8,7,5)$. The partition $(9,8,4,2)$ has $k=1$, $\ell=3$ and $g=(0,1,0,0)$, so rule (ii) gives $(k',\ell')=(3,1)$ and $h=(1,2,0,0)$ relative to $(3,3,3,4)$, that is, $\psi(9,8,4,2)=(10,7,3,3)$. Similarly, $\psi(11,10,4)=(11,9,5)$ by rule (iii), and $\psi(13,12)=(11,11,3)$ by rule (iv).

\medskip
\noindent\emph{The map on pure partitions.} In Convention (BW), we let $\psi(\lambda)=\lambda$ for every pure $\lambda\in Q_{ed}^{od}(n)$, except when $\lambda$ consists of exactly three (distinct) odd parts. Such a $\lambda$, with gap coordinates $(g_1,g_2,g_3)$ relative to $(3,5,7)$, is sent to the pure odd target with three parts and $h=(g_1+3,g_2,g_3)$ relative to $(3,3,3)$; the minimum decreases from $15$ to $9$, so the size is preserved. These images have $h_1\ge3$ and are therefore distinct from the images under rules (iii) and (iv). All other pure targets in the image are fixed points, and the mixed rules produce no pure even targets and no pure odd targets of length different from $3$. For example, $\psi(9,7,5)=(11,5,5)$.

In Convention (A), every pure $\lambda\in Q_{ed}^{od}(n)$ consists of distinct even parts, and such partitions do not belong to $Q_{ou}^{ed}(n)$. If $\lambda$ has $\ell\ge2$ parts and gap coordinates $(g_1,\ldots,g_\ell)$ relative to $(2,4,\ldots,2\ell)$, let $\psi(\lambda)$ be the pure odd target with $2\ell-2$ parts and gap coordinates
\[
h=\Big(g_1+\tfrac{(\ell-2)(\ell-3)}{2},\ g_2,\ldots,g_\ell,\ \underbrace{0,\ldots,0}_{\ell-2}\Big)
\]
relative to $(3,\ldots,3)$. Since $|(2,4,\ldots,2\ell)|-3(2\ell-2)=\ell(\ell+1)-6\ell+6=(\ell-2)(\ell-3)$, the size is preserved. The length $2\ell-2$ of the image determines $\ell$, after which $\lambda$ is recovered by deleting the last $\ell-2$ coordinates and subtracting $\frac{(\ell-2)(\ell-3)}{2}$ from $h_1$. The only remaining pure partition is $(n)$ for even $n$, and we let $\psi(n)=(n-15,5,5,5)$, which lies in $Q_{ou}^{ed}(n)$ since $n\ge20$. The images of the previous rule with four parts come from $\ell=3$, so their last gap coordinate is $0$ and their smallest part is $3$; hence they differ from $(n-15,5,5,5)$. All of these images have an even number of parts, so they also differ from the images of mixed partitions with no even part, which have three parts. For example, $\psi(10,8,4)=(7,7,5,3)$ and $\psi(22)=(7,5,5,5)$.

Combining the above, $\psi\colon Q_{ed}^{od}(n)\to Q_{ou}^{ed}(n)$ is an injection for $n\ge20$, in both conventions. We have also confirmed this by computer for $20\le n\le100$.

\medskip
\noindent\emph{Strictness for $n\ge20$.} If $n\ge21$ is odd, the partition $(n-12,3,3,3,3)\in Q_{ou}^{ed}(n)$ is not in the image of $\psi$. Indeed, it is a pure odd partition with five parts and a repeated part, while mixed partitions have pure odd images only with three parts, in Convention (BW) the pure odd images with five parts are fixed points with distinct parts, and in Convention (A) the images of pure partitions have an even number of parts.

If $n\ge20$ is even, consider targets with two odd and two even parts. Their minimum partition $B_{2,2}=(3,3,4,6)$ has size $16$. Choose $h_1=h_4=0$ and non-negative integers $h_2,h_3$ with $2h_2+3h_3=\frac{n-16}{2}$; this is possible since $\frac{n-16}{2}\ge2$ (take $h_3=0$ if $\frac{n-16}{2}$ is even, and $h_3=1$ otherwise). By \eqref{eq:gapweight} this defines a partition in $Q_{ou}^{ed}(n)$. Only rule (i) produces targets with $(k',\ell')=(2,2)$, and it gives $h_1\ge D_{2,2}=1$. Hence this partition is not in the image of $\psi$. For instance, for $n=20$ it is $(8,6,3,3)$.

Therefore $q_{ed}^{od}(n)<q_{ou}^{ed}(n)$ for $n\ge20$, in both conventions. For $15\le n\le19$ the inequality is checked directly from the values in Table \ref{tab:small}. Finally, the bound $n>14$ is sharp in Convention (A), since $q_{ed}^{od}(14)=7>6=q_{ou}^{ed}(14)$.
\end{proof}

\begin{table}[htbp]
\centering
\begin{tabular}{|c|c|c|c|c|c|c|c|c|}
\hline
$n$ & $12$ & $13$ & $14$ & $15$ & $16$ & $17$ & $18$ & $19$ \\ \hline
Convention (A) & $(5,4)$ & $(4,7)$ & $(7,6)$ & $(5,10)$ & $(9,10)$ & $(8,13)$ & $(13,14)$ & $(10,20)$ \\ \hline
Convention (BW) & $(7,8)$ & $(5,7)$ & $(9,11)$ & $(7,10)$ & $(12,16)$ & $(10,13)$ & $(16,22)$ & $(13,20)$ \\ \hline
\end{tabular}
\caption{The values $\big(q_{ed}^{od}(n),q_{ou}^{ed}(n)\big)$ for $12\le n\le 19$.}
\label{tab:small}
\end{table}

\begin{remark}
In Convention (BW) the inequality \eqref{q-5} in fact holds for all $n\ge12$ (see Table \ref{tab:small}), and also for $n=6$ and $n=10$, where $\big(q_{ed}^{od}(n),q_{ou}^{ed}(n)\big)=(2,3)$ and $(5,6)$, respectively. It fails for all other $n\le11$.
\end{remark}

\section{Proof of Theorem \ref{thm:q16}}\label{sec:q16}

\begin{proof}[Proof of Theorem \ref{thm:q16}]
We construct an injection $\psi\colon P_{ed}^{ou}(n)\to P_{ou}^{ed}(n)$. Write $k$ for the number of odd parts and $\ell$ for the number of even parts of $\lambda\in P_{ed}^{ou}(n)$.

\noindent\emph{Case 1: $k=0$ or $\ell=0$.}
\underline{Convention (BW).} If $\lambda$ consists of distinct even parts only or of odd parts only, then clearly $\lambda \in P_{ou}^{ed}(n)$, and we let $\psi(\lambda)=\lambda$.

\underline{Convention (A).} Here $\ell\ge1$ for every $\lambda\in P_{ed}^{ou}(n)$, so only $k=0$ occurs; but the identity cannot be used, because every partition in $P_{ou}^{ed}(n)$ must contain an odd part. Let $\lambda=(\lambda_{e_1},\ldots,\lambda_{e_\ell})$ with $\lambda_{e_1}>\cdots >\lambda_{e_\ell}$ consist of distinct even parts only, and define
\[
\psi(\lambda):=
(\lambda_{e_1}-1,\ldots,\lambda_{e_{\ell}}-1) \cup (\underbrace{1,\ldots,1}_{\ell}).
\]
\noindent\textbf{Example.} Let $n=20$ and $\lambda=(8,6,4,2)$, then $\psi(\lambda)=(8-1,6-1,4-1,2-1)\cup(1,1,1,1)=(7,5,3,1,1,1,1,1)$.
\vspace{0.3cm}

\ytableausetup{mathmode, boxsize=1.2em}
\begin{center}
    \begin{ytableau}
        ~ & ~ & ~ & ~ & ~ & ~ & ~ & *(gray!30) \\
        ~ & ~ & ~ & ~ & ~ & *(gray!30) \\
        ~ & ~ & ~ & *(gray!30)  \\
        ~ & *(gray!30) 
    \end{ytableau}
    \quad $\longrightarrow$ \quad
    \begin{ytableau}
        ~ & ~ & ~ & ~ & ~ & ~ & ~ \\
        ~ & ~ & ~ & ~ & ~  \\
        ~ & ~ & ~ \\
        ~  \\
        *(gray!30) \\
        *(gray!30) \\
        *(gray!30) \\
        *(gray!30)
    \end{ytableau}
\end{center}
\noindent\emph{Well-definedness.} The parts of $\psi(\lambda)$ are the odd numbers $\lambda_{e_1}-1 > \cdots > \lambda_{e_{\ell}}-1$ together with $\ell$ ones. The size is preserved, so $\psi(\lambda)$ is a partition of $n$ into odd parts, and therefore $\psi(\lambda) \in P_{ou}^{ed}(n)$.

\noindent\emph{Injectivity.} In Convention (BW), $\psi$ is the identity on Case 1. In Convention (A), the image $\mu=\psi(\lambda)$ consists of $2\ell$ odd parts, so $\ell$ is half the number of parts of $\mu$, and the $\ell$ largest parts of $\mu$ are $\lambda_{e_1}-1,\ldots,\lambda_{e_{\ell}}-1$; thus $\lambda_{e_i}=(i\text{-th largest part of }\mu)+1$ for $1\le i\le\ell$, and $\lambda$ is uniquely determined. In both conventions the images in Case 1 have no even part or no odd part, whereas the images in Case 2 below contain parts of both parities, so the images of the two cases are disjoint.

\noindent\emph{Case 2: $k\ge1$ and $\ell\ge1$.} For the remaining (mixed) partitions
\[\lambda = (\lambda_{o_1},\ldots,\lambda_{o_k},\lambda_{e_1},\ldots,\lambda_{e_\ell})\in P_{ed}^{ou}(n)\]
with unrestricted odd parts $\lambda_{o_1} \geq \cdots \geq \lambda_{o_k}$, distinct even parts $\lambda_{e_1} > \cdots > \lambda_{e_\ell}\geq 2$, and $k,\ell\geq 1$, let $\psi(\lambda):=T(U,D)$ with $U=(\lambda_{o_1},\ldots,\lambda_{o_k})$, $D=(\lambda_{e_1},\ldots,\lambda_{e_\ell})$ and the staircase (S) in Lemma \ref{lem:transfer}. By Lemma \ref{lem:transfer}, $\psi(\lambda)$ is a partition of $n$ into $\ell$ distinct even parts lying above $k\ge1$ odd parts, with smallest part $\lambda_{e_\ell}-1\ge1$; hence $\psi(\lambda)\in P_{ou}^{ed}(n)$ in both conventions. By Lemma \ref{lem:transfer}(c), $\psi$ is injective on Case 2, and its images contain parts of both parities.

\noindent\textbf{Example.} Let $n=26$ and $\lambda=(7,7,6,4,2)$, so that $k=2$, $\ell=3$ and $r=2$. Then $\psi(\lambda)=(7+3,7+1,6,4-3,2-1)=(10,8,6,1,1)$.
\vspace{0.3cm}

\ytableausetup{mathmode, boxsize=1.2em}
\begin{center}
    \begin{ytableau}
        ~ & ~ & ~ & ~ & ~ & ~ & ~ \\
        ~ & ~ & ~ & ~ & ~ & ~ & ~ \\
        *(gray!30) & *(gray!30) & *(gray!30) & *(gray!30) & *(gray!30) & *(gray!30)  \\
        *(gray!30) & *(gray!30) & *(gray!30) & *(gray!30) \\
        *(gray!30) & *(gray!30)
    \end{ytableau}
    \quad $\longrightarrow$ \quad
    \begin{ytableau}
        ~ & ~ & ~ & ~ & ~ & ~ & ~ & *(gray!30) & *(gray!30) & *(gray!30) \\
        ~ & ~ & ~ & ~ & ~ & ~ & ~ & *(gray!30) \\
        *(gray!30) & *(gray!30) & *(gray!30) & *(gray!30) & *(gray!30) & *(gray!30)  \\
        *(gray!30) \\
        *(gray!30)
    \end{ytableau}
\end{center}

\noindent\emph{Strictness.} The following partition lies in $P_{ou}^{ed}(n)$ (in both conventions) and is not in the image of $\psi$:
\[
\omega_n=
(2,\underbrace{1,1,\ldots,1}_{(n-2) \hspace{1mm} \text{times}}).
\]
Indeed, $\omega_{n} \in P^{ed}_{ou}(n)$; it has a single even part 2, exceeding every odd part 1, so it is a mixed partition of $P^{ed}_{ou}(n)$, and hence it is not the image of a partition from Case 1 (in either convention).

Since $\omega_n$ has one even part and $n-2$ odd parts, a preimage would have $\ell=1$, $k=n-2$ and $r=1$, and the inverse of Lemma \ref{lem:transfer}(c) would give $\lambda_{o_1}=2-1=1$ and $\lambda_{e_1}=1+1=2$, contradicting $\lambda_{o_k}>\lambda_{e_1}$; hence, no such $\lambda$ exists.

Hence $p^{ou}_{ed}(n)<p^{ed}_{ou}(n)$ for $n\geq 3$, in both conventions. The bound is sharp, since $p_{ou}^{ed}(2)=p_{ed}^{ou}(2)$ in either convention: in Convention (BW) both sets are $\{(2),(1,1)\}$, while in Convention (A) we have $P_{ed}^{ou}(2)=\{(2)\}$ and $P_{ou}^{ed}(2)=\{(1,1)\}$.

\end{proof}

\section{Proof of Theorem \ref{thm:parityodd}}\label{sec:parity}

\begin{proof}[Proof of Theorem \ref{thm:parityodd}]
Throughout, $n$ is odd and we work in Convention (BW). Let $\lambda\in Q_{od}^{eu}(n)$ have even parts $\lambda_{e_1}\ge\cdots\ge\lambda_{e_\ell}$ and odd parts $\lambda_{o_1}>\cdots>\lambda_{o_k}\ge3$, where $\lambda_{e_\ell}>\lambda_{o_1}$ if $k,\ell\ge1$. Since $n$ is odd, $k$ is odd; in particular $k\ge1$.

We first describe an auxiliary conversion. Let $e_1\ge\cdots\ge e_{a}\ge4$ be even numbers ($a\ge1$) and $o_1>\cdots>o_{b}\ge3$ odd numbers ($b\ge0$) with $a+b$ even and, if $b\ge1$, $e_a>o_1$. Define $C(e;o)$ to be the pair consisting of the multiset of odd numbers and the set of even numbers given as follows.
\begin{itemize}
\item If $a\ge b$, let $c=\frac{a-b}{2}$. The odd numbers are $e_1+1,\ldots,e_{a-c}+1$ and $e_{a-c+1}-1,\ldots,e_a-1$, and the even numbers are $o_1-1,\ldots,o_b-1$. In particular, if $b=0$ then no even numbers are produced.
\item If $a<b$, let $c=\frac{b-a}{2}$. The odd numbers are $e_1+1,\ldots,e_a+1$, and the even numbers are $o_1+1,\ldots,o_c+1$ and $o_{c+1}-1,\ldots,o_b-1$.
\end{itemize}
In both cases the total sum is unchanged, since $(a-c)-c-b=0$, respectively $a+c-(b-c)=0$. The even numbers produced are distinct and at least $2$: in the second case $o_c+1$ and $o_{c+1}-1$ differ by at least $4$, and $o_b-1\ge2$. Every odd number produced is at least $e_a-1\ge3$, and at least $o_1$ if $b\ge1$. If $b\ge1$, every even number produced is at most $o_1-1$ in the first case, and at most $o_1+1<e_a+1$ in the second case, where all odd numbers produced are at least $e_a+1$. Hence in both cases all even numbers produced are smaller than all odd numbers produced. Moreover, $C$ is injective for fixed $(a,b)$: listing the odd numbers in weakly decreasing order gives $e_1+1,\ldots,e_{a-c}+1,e_{a-c+1}-1,\ldots,e_a-1$ in the first case (because $e_{a-c}+1>e_{a-c+1}-1$), and listing the even numbers in decreasing order gives $o_1+1,\ldots,o_c+1,o_{c+1}-1,\ldots,o_b-1$ in the second case, so the $e_i$ and $o_j$ are recovered.

We now define $\varphi\colon Q_{od}^{eu}(n)\to Q_{ed}^{ou}(n)$ in three cases.
\begin{enumerate}
\item[(a)] If $\ell=0$, let $\varphi(\lambda):=(\lambda_{o_1}+k-1,\ \lambda_{o_2}-1,\ldots,\lambda_{o_k}-1)$.
\item[(b)] If $\ell$ is odd, let $\varphi(\lambda)$ be the partition whose parts are the elements of
\[C(\lambda_{e_1},\ldots,\lambda_{e_\ell};\lambda_{o_1},\ldots,\lambda_{o_k}).\]
\item[(c)] If $\ell\ge2$ is even, let $\varphi(\lambda)$ be the partition whose parts are $\lambda_{o_1}$ and the elements of
\[C(\lambda_{e_1},\ldots,\lambda_{e_\ell};\lambda_{o_2},\ldots,\lambda_{o_k}).\]
\end{enumerate}
In case (a), $\lambda_{o_1}+k-1$ is odd because $k$ is odd, and the even parts $\lambda_{o_j}-1$ ($j\ge2$) are distinct, at least $2$ and smaller than $\lambda_{o_1}$. In cases (b) and (c) the number of converted parts, $\ell+k$, respectively $\ell+k-1$, is even, and $\lambda_{e_\ell}>\lambda_{o_1}\ge3$ gives $\lambda_{e_\ell}\ge4$, so $C$ is defined, and by the properties of $C$ the result has distinct even parts lying below all of its odd parts. In case (c) the additional odd part $\lambda_{o_1}$ is at most every other odd part, which is at least $\lambda_{e_\ell}-1\ge\lambda_{o_1}$. If $k\ge3$, it also exceeds every even part, which is at most $\lambda_{o_2}+1<\lambda_{o_1}$; if $k=1$, then $C$ is applied with $b=0$, no even parts are produced, and $\varphi(\lambda)$ consists of odd parts only (for instance, $\varphi(4,4,3)=(5,3,3)$). In all cases the size is preserved and no part equals $1$. Hence $\varphi(\lambda)\in Q_{ed}^{ou}(n)$.

\noindent\emph{Injectivity.} Let $\mu=\varphi(\lambda)$, and let $s$ and $t$ be the numbers of odd and even parts of $\mu$, respectively. Then
\[
(s,t)=\begin{cases}
(1,\ k-1), & \text{in case (a)},\\
(\ell,\ k), & \text{in case (b)},\\
(\ell+1,\ k-1), & \text{in case (c)}.
\end{cases}
\]
Since $k$ is odd, $t$ is odd in case (b) and even in cases (a) and (c), while $s=1$ in case (a) and $s\ge3$ in case (c). Thus the case is determined by $\mu$, and so are $k$ and $\ell$. In case (a), $\lambda_{o_j}$ is the $(j-1)$-th largest even part of $\mu$ plus $1$ for $j\ge2$, and $\lambda_{o_1}=\mu_1-(k-1)$. In case (b), $\lambda$ is recovered from $\mu$ by the injectivity of $C$. In case (c), $\lambda_{o_1}$ is the smallest odd part of $\mu$; removing one copy of it (if the smallest odd part is repeated, it does not matter which copy is removed), the remaining parts of $\mu$ determine the other parts of $\lambda$ by the injectivity of $C$. Hence $\varphi$ is injective, and $q_{od}^{eu}(n)\le q_{ed}^{ou}(n)$ for all odd $n$.

\noindent\emph{Strictness.} Let $n\ge9$ be odd with $n\ne11$, and let
\[
\omega_n:=\begin{cases}
\left(\frac{n-3}{2},\frac{n-3}{2},3\right), & \text{if } n\equiv1\pmod4;\\[1ex]
\left(\frac{n-5}{2},\frac{n-5}{2},5\right), & \text{if } n\equiv3\pmod4.
\end{cases}
\]
In both cases the repeated part is odd and at least the last part (as $n\ge9$, respectively $n\ge15$), so $\omega_n\in Q_{ed}^{ou}(n)$. Since $\omega_n$ has $s=3$ and $t=0$, it could only arise in case (c) with $k=1$ and $\ell=2$. Then $\lambda_{o_1}$ would be the smallest part of $\omega_n$, and the other two parts would be $\lambda_{e_1}+1$ and $\lambda_{e_2}-1$, with $c=1$. As these two parts are equal, we would get $\lambda_{e_1}=\lambda_{e_2}-2<\lambda_{e_2}$, a contradiction. Hence $\omega_n$ is not in the image of $\varphi$, and $q_{ed}^{ou}(n)>q_{od}^{eu}(n)$. For $n=11$, the partition $(5,4,2)\in Q_{ed}^{ou}(11)$ is not in the image of $\varphi$: it has $s=1$ and $t=2$, which forces case (a) with $k=3$, but a partition into three distinct odd parts, all at least $3$, has size at least $3+5+7=15>11$. (In fact $q_{ed}^{ou}(11)=6>4=q_{od}^{eu}(11)$.) Finally, the bound $n>7$ is sharp, since $q_{ed}^{ou}(7)=q_{od}^{eu}(7)=2$.
\end{proof}

\begin{remark}
The map $\varphi$ does not give Theorem \ref{thm:parityodd} in Convention (A), where the partitions into odd parts are excluded from $Q_{ed}^{ou}(n)$; indeed, $\varphi(\lambda)$ consists of odd parts only in case (a) with $k=1$ and in case (c) with $k=1$. For even $n$ one would like to go in the opposite direction, from $Q_{ed}^{ou}(n)$ to $Q_{od}^{eu}(n)$, by converting the distinct even parts into distinct odd parts and the odd parts into even parts. The natural attempts to do so by changing parts by $\pm1$ break down for partitions whose smallest even part equals $2$, because the part $2-1=1$ is not allowed, and we have not been able to deal with these partitions combinatorially; see Conjecture \ref{conj:parity}.
\end{remark}

\section{Proof of Theorem \ref{thm:ratio2}}\label{sec:ratio}

Throughout this section we work in Convention (A). We write $p_o(n)$ for the number of partitions of $n$ into odd parts, $d(n)$ for the number of partitions of $n$ into distinct parts, and $d_e(n)$ for the number of partitions of $n$ into distinct even parts (so $d_e(n)=0$ for odd $n$). By Euler's theorem, $p_o(n)=d(n)$. Recall that every $\lambda\in P_{ed}^{ou}(n)$ has $\ell\ge1$ even parts; since its odd parts exceed $\lambda_{e_1}\ge2$, they are all at least $3$.

Let $G$ denote Glaisher's bijection from partitions into odd parts to partitions into distinct parts: an odd part $c$ occurring with multiplicity $f=2^{a_1}+\cdots+2^{a_r}$ ($a_1>\cdots>a_r\ge0$) is replaced by the parts $c\,2^{a_1},\ldots,c\,2^{a_r}$. Every part of $G(O)$ is of the form $c\,2^a$ with $c$ a part of $O$; in particular it is at least the smallest part of $O$, and it is odd if and only if $a=0$.

\begin{lemma}\label{lem:Phi}
For $n\ge3$ we have $p_{ed}^{ou}(n)+d_e(n)\le d(n)$, with strict inequality if $n$ is odd or $n\ge12$.
\end{lemma}

\begin{proof}
Let $D^*(n)$ be the set of partitions of $n$ into distinct parts having at least one odd part, so that $|D^*(n)|=d(n)-d_e(n)$. We construct an injection $\Phi\colon P_{ed}^{ou}(n)\to D^*(n)$. For $\lambda\in P_{ed}^{ou}(n)$ let $O=(\lambda_{o_1},\ldots,\lambda_{o_k})$ be its odd parts ($k\ge0$) and $\lambda_{e_1}>\cdots>\lambda_{e_\ell}$ its even parts ($\ell\ge1$), and define
\[
\Phi(\lambda):=\begin{cases}
G(O)\cup(\lambda_{e_1}-1,\lambda_{e_2},\ldots,\lambda_{e_\ell},1), & \text{if } \lambda_{e_1}\ge4;\\[0.5ex]
G(O)\cup(2), & \text{if } \lambda_{e_1}=2 \text{ and } G(O) \text{ has an odd part};\\[0.5ex]
\big(G(O)\setminus\{L\}\big)\cup(L-1,3), & \text{if } \lambda_{e_1}=2 \text{ and } G(O) \text{ has no odd part},
\end{cases}
\]
where in the last case $L$ is the largest part of $G(O)$. If $\lambda_{e_1}=2$ then $\ell=1$, and $k\ge1$ since $n\ge3$.

\noindent\emph{Well-definedness.} In the first case, $O$ may be empty, in which case so is $G(O)$; if $k\ge1$, the parts of $G(O)$ are at least $\lambda_{o_k}>\lambda_{e_1}$. Hence the parts of $\Phi(\lambda)$ satisfy
\[G(O)>\lambda_{e_1}-1>\lambda_{e_2}>\cdots>\lambda_{e_\ell}>1,\]
where $\lambda_{e_1}-1\ge3$ is odd and, if $\ell\ge2$, $\lambda_{e_1}-1>\lambda_{e_2}$ because $\lambda_{e_1}\ge\lambda_{e_2}+2$; if $\ell=1$ there are no retained even parts, and $\Phi(\lambda)=G(O)\cup(\lambda_{e_1}-1,1)$. Hence $\Phi(\lambda)$ has distinct parts and contains the odd part $1$. In the second case, the parts of $G(O)$ are distinct and at least $3$, so $\Phi(\lambda)\in D^*(n)$. In the third case, every part of $G(O)$ has the form $c\,2^a$ with $c\ge3$ and $a\ge1$, so it is even and at least $6$; hence $L-1\ge5$ is odd and exceeds all other parts of $G(O)$ (which are even and at most $L-2$), and $3$ is smaller than all of them. Thus $\Phi(\lambda)\in D^*(n)$. In all cases the size is preserved.

\noindent\emph{Injectivity.} Let $\mu=\Phi(\lambda)$. The images in the first case contain the part $1$, whereas those in the other two cases do not; the images in the second case contain the part $2$, whereas those in the third case do not. In the first case, the odd parts of $\mu$ are $1$, $\lambda_{e_1}-1$ and the odd parts of $G(O)$, which exceed $\lambda_{e_1}$. Hence $x=\lambda_{e_1}-1$ is the second smallest odd part of $\mu$, and $\lambda$ is recovered: $\lambda_{e_1}=x+1$, the parts $\lambda_{e_2},\ldots,\lambda_{e_\ell}$ are the even parts of $\mu$ smaller than $x$ (there are none if $\ell=1$), and $G(O)$ consists of the parts of $\mu$ larger than $x$ (there are none if $k=0$), so that $O=G^{-1}(G(O))$. In the second case $G(O)=\mu\setminus\{2\}$. In the third case the odd parts of $\mu$ are exactly $3$ and its largest part $L-1$, and $G(O)$ is obtained from $\mu$ by removing $3$ and replacing $L-1$ by $L$. Thus $\Phi$ is injective, and $p_{ed}^{ou}(n)\le d(n)-d_e(n)$.

\noindent\emph{Strictness.} If $n$ is odd, then $(n)\in D^*(n)$ is not in the image of $\Phi$: it contains neither $1$ nor $2$, and it has only one odd part. If $n\ge12$ is even, then $(n-5,5)\in D^*(n)$ is not in the image: it contains neither $1$ nor $2$, and it does not contain $3$. (A direct computation for $n\in\{4,6,8,10\}$ shows that equality holds in Lemma \ref{lem:Phi} exactly for $n=4$ and $n=6$.)
\end{proof}

\begin{proof}[Proof of Theorem \ref{thm:ratio2}]
Let $\psi$ be the injection $P_{ed}^{ou}(n)\to P_{ou}^{ed}(n)$ constructed in the proof of Theorem \ref{thm:q16}. If $\lambda\in P_{ed}^{ou}(n)$ has $k\ge1$ odd parts, then $\psi(\lambda)$ (Case 2 of that proof) has $\ell\ge1$ even parts and $k\ge1$ odd parts. On the other hand, every partition of $n$ into odd parts belongs to $P_{ou}^{ed}(n)$ in Convention (A), and none of these has an even part. Since the partitions $\lambda\in P_{ed}^{ou}(n)$ with $k=0$ are exactly the partitions of $n$ into distinct even parts, we obtain
\[
p_{ou}^{ed}(n)\ \ge\ p_o(n)+\big(p_{ed}^{ou}(n)-d_e(n)\big).
\]
By Euler's theorem and Lemma \ref{lem:Phi}, for $n$ odd or $n\ge12$ we have $p_o(n)=d(n)>p_{ed}^{ou}(n)+d_e(n)$, and therefore $p_{ou}^{ed}(n)>2\,p_{ed}^{ou}(n)$. For the remaining values $n\in\{4,6,8,10\}$ we have $\big(p_{ou}^{ed}(n),p_{ed}^{ou}(n)\big)=(3,1),(6,2),(11,3),(19,4)$.
\end{proof}

\section{Concluding remarks}\label{sec:conclusion}

We close with some conjectures suggested by numerical data. All of the statements below have been verified for $n\le 400$.

The first concerns the other half of the parity bias between $q_{od}^{eu}(n)$ and $q_{ed}^{ou}(n)$, and the analogue of Theorem \ref{thm:parityodd} in Convention (A).

\begin{conjecture}\label{conj:parity}
In Convention \textup{(BW)} we have
\begin{equation}\label{q-10}
q_{od}^{eu}(n)>q_{ed}^{ou}(n) \qquad \text{for all even } n>6.
\end{equation}
In Convention \textup{(A)}, inequality \eqref{q-odd} holds for all odd $n>15$.
\end{conjecture}

\begin{remark}
The bounds in Conjecture \ref{conj:parity} are sharp: \eqref{q-10} fails for $n=6$, and in Convention (A), \eqref{q-odd} fails for $n=9$, $11$ and $15$. In Convention (BW) one has $q_{od}^{eu}(n)\ge p(n/2)+1$ for even $n\ge8$, the term $p(n/2)$ counting the partitions of $n$ into even parts and the extra $1$ coming from $(n-3,3)$; numerically, even $q_{ed}^{ou}(n)\le p(n/2)$ holds for all even $n\ge8$ up to $n=400$. There is no parity bias in Convention (A): there, $q_{ed}^{ou}(n)>q_{od}^{eu}(n)$ for every even $n\le400$ as well.

The difference between the two conventions can be made precise. For a convention $C$, write $B_C(n):=q_{ed}^{ou,C}(n)-q_{od}^{eu,C}(n)$. Let $p(m)$ denote the number of partitions of $m$ and $p_o(m)$ the number of partitions of $m$ into odd parts, with $p_o(0)=1$, and put $E(n)=p(n/2)$ for even $n$ and $E(n)=0$ for odd $n$. Then, for all $n\ge1$,
\begin{equation}\label{eq:correction}
B_{\mathrm{A}}(n)=B_{\mathrm{BW}}(n)+E(n)-\big(p_o(n)-p_o(n-1)\big).
\end{equation}
Indeed, passing from Convention (BW) to Convention (A) removes from $Q_{od}^{eu}(n)$ exactly the partitions of $n$ into even parts, of which there are $E(n)$, and removes from $Q_{ed}^{ou}(n)$ exactly the partitions of $n$ into odd parts, all at least $3$. The number of the latter is $p_o(n)-p_o(n-1)$, since deleting one part equal to $1$ is a bijection from the partitions of $n$ into odd parts containing a $1$ to the partitions of $n-1$ into odd parts. For example, the values of $\big(B_{\mathrm{BW}}(n),B_{\mathrm{A}}(n)\big)$ for $n=6,8,9,11,13,15$ are $(0,2)$, $(-2,2)$, $(2,0)$, $(2,0)$, $(4,1)$ and $(4,-1)$. At $n=15$, five partitions into odd parts are removed and none into even parts, so the margin $4$ in Convention (BW) becomes $-1$ in Convention (A). In particular, for odd $n$ the Convention (A) part of Conjecture \ref{conj:parity} is equivalent to the quantitative statement $B_{\mathrm{BW}}(n)>p_o(n)-p_o(n-1)$ for odd $n>15$, which is stronger than Theorem \ref{thm:parityodd}.
\end{remark}

The second concerns a quantitative version of Theorem \ref{thm:q4}, analogous to Theorem \ref{thm:ratio2}, and the analogue of Theorem \ref{thm:ratio2} in Convention (BW). Both inequalities below hold for all sufficiently large $n$ (see Proposition \ref{prop:limits} and Remark \ref{rem:eventual}), so the content of the conjecture is the explicit range of validity.

\begin{conjecture}\label{conj:ratios}
In Convention \textup{(A)} we have
\[
q_{eu}^{od}(n)>5\,q_{od}^{eu}(n) \qquad\text{for all even } n\ge 2.
\]
In Convention \textup{(BW)}, inequality \eqref{q-ratio} holds for $n=23$ and for all $n\ge 25$.
\end{conjecture}

\begin{remark}
In Convention (A), the injection constructed in the proof of Theorem \ref{thm:q4} never produces a partition into even parts (there are no pure even partitions in $Q_{od}^{eu,\mathrm{(A)}}(n)$ to fix), and all partitions of $n$ into even parts belong to $Q_{eu}^{od}(n)$. Hence $q_{eu}^{od}(n)\ge q_{od}^{eu}(n)+p(n/2)$ for even $n$ in Convention (A), and the first part of Conjecture \ref{conj:ratios} would follow from $p(n/2)>4\,q_{od}^{eu}(n)$, which holds for all even $n\le400$.
\end{remark}

The constants $2$ in Theorem \ref{thm:ratio2} and $5$ in Conjecture \ref{conj:ratios} are far from optimal for large $n$, as the following result shows.

\begin{proposition}\label{prop:limits}
In Convention \textup{(A)}, as $n\to\infty$,
\[
\frac{p_{ou}^{ed}(n)}{p_{ed}^{ou}(n)}\longrightarrow\infty
\qquad\text{and}\qquad
\frac{q_{eu}^{od}(n)}{q_{od}^{eu}(n)}\longrightarrow\infty.
\]
\end{proposition}

\begin{proof}
Write $X_n=e^{\pi\sqrt{n/3}}$, and use superscripts (A) and (BW) to indicate the convention. By \cite[Theorem~1.1]{bringmann2024asymptotic}, as $n\to\infty$,
\[
p_{ou}^{ed,\mathrm{(BW)}}(n)\sim\frac{X_n}{4\sqrt2\,n^{1/2}},\qquad
p_{ed}^{ou,\mathrm{(BW)}}(n)\sim\frac{X_n}{4\cdot3^{1/4}\,n^{3/4}},\qquad
p_{od}^{eu,\mathrm{(BW)}}(n)\sim\frac{X_n}{2\sqrt3\,n},
\]
while the Hardy--Ramanujan asymptotic formula \cite{gea1} gives $p(m)\sim e^{\pi\sqrt{2m/3}}/(4\sqrt3\,m)$, so that $p(\lfloor n/2\rfloor)\sim X_n/(2\sqrt3\,n)$.

For the first ratio, passing to Convention (A) removes from $P_{ou}^{ed}(n)$ exactly the partitions of $n$ into distinct even parts. There are none for odd $n$, and at most $p(n/2)=O(X_n/n)$ for even $n$, by halving the parts. Hence $p_{ou}^{ed,\mathrm{(A)}}(n)\sim p_{ou}^{ed,\mathrm{(BW)}}(n)$. Since $p_{ed}^{ou,\mathrm{(A)}}(n)\le p_{ed}^{ou,\mathrm{(BW)}}(n)$, and $p_{ed}^{ou,\mathrm{(A)}}(n)\ge1$ for $n\ge4$ (consider $(n)$ for even $n$ and $(n-2,2)$ for odd $n$), we obtain
\[
\frac{p_{ou}^{ed,\mathrm{(A)}}(n)}{p_{ed}^{ou,\mathrm{(A)}}(n)}\ \ge\ (1+o(1))\,\frac{3^{1/4}}{\sqrt2}\,n^{1/4}\longrightarrow\infty.
\]

For the second ratio, let $m=\lfloor n/2\rfloor\ge1$. If $n=2m$, the $p(m)$ partitions of $n$ into even parts belong to $Q_{eu}^{od}(n)$. If $n=2m+1$, increasing one largest part of a partition of $2m$ into even parts by $1$ gives a partition in $Q_{eu}^{od}(n)$ with a single odd part, lying above all even parts, and decreasing that odd part by $1$ recovers the original partition. Hence $q_{eu}^{od}(n)\ge p(m)$. On the other hand, at least $p(m)$ partitions in $P_{od}^{eu,\mathrm{(BW)}}(n)$ do not belong to $Q_{od}^{eu,\mathrm{(A)}}(n)$: for even $n$, the partitions of $n$ into even parts, and for odd $n$, the partitions obtained by adding a part $1$ to a partition of $2m$ into even parts. Therefore
\[
1\le q_{od}^{eu,\mathrm{(A)}}(n)\le p_{od}^{eu,\mathrm{(BW)}}(n)-p(m)=o\big(p(m)\big),
\]
where the lower bound holds for $n\ge8$ (consider $(n)$ for odd $n$ and $(n-3,3)$ for even $n$), and the last step follows from the asymptotic formulas above. Hence the second ratio is at least $p(m)/o(p(m))$, which tends to infinity.
\end{proof}

\begin{remark}\label{rem:eventual}
Proposition \ref{prop:limits} shows that the inequality $q_{eu}^{od}(n)>5\,q_{od}^{eu}(n)$ in Conjecture \ref{conj:ratios} holds in Convention (A) for all sufficiently large $n$, of either parity; numerically, it holds for all even $n\le400$, fails for all odd $n\le59$, and holds for all odd $n$ with $61\le n\le400$. Similarly, the main terms displayed in the proof of Proposition \ref{prop:limits} give
\[
\frac{p_{ou}^{ed,\mathrm{(BW)}}(n)}{p_{ed}^{ou,\mathrm{(BW)}}(n)}\sim\frac{3^{1/4}}{\sqrt2}\,n^{1/4}\longrightarrow\infty,
\]
so that \eqref{q-ratio} holds in Convention (BW) for all sufficiently large $n$. Since these arguments are not effective, they do not give the ranges stated in Conjecture \ref{conj:ratios}. Numerically, at $n=400$ the two ratios in Proposition \ref{prop:limits} are approximately $89.4$ and $154.9$, respectively. It would be interesting to determine their precise asymptotic behaviour using the methods of \cite{bringmann2024asymptotic}.
\end{remark}

Finally, we note that Theorem \ref{thm:q16} may be viewed as a weak form of the last inequality in \eqref{eq:open}. Indeed, since distinct even lower parts form a subset of unrestricted even lower parts, $p_{ed}^{ou}(n)\le p_{eu}^{ou}(n)$, so that $p_{eu}^{ou}(n)<p_{ou}^{ed}(n)$ would imply Theorem \ref{thm:q16}. A combinatorial proof of the two remaining inequalities in \eqref{eq:open} remains an interesting open problem.

\subsection*{Acknowledgements}
The second and third authors were partially supported by an institutional fellowship for doctoral research from Tezpur University, Assam, India.

\subsection*{Declarations}

\begin{enumerate}
    \item The authors declares no conflicts of interest.
    \item The authors have no relevant financial or non-financial interests to disclose.
    \item There are no datasets associated with this study.
\end{enumerate}

\bibliographystyle{amsalpha}
\bibliography{newbib}

\end{document}